\documentclass[10pt]{amsart}

\usepackage{amssymb,amsthm,amsthm,amsfonts,latexsym}
\usepackage{amsmath}
\usepackage{mathrsfs}
\usepackage{stmaryrd}
\usepackage{accents}
\usepackage[utf8]{inputenc}
\usepackage[]{hyperref}

\usepackage[paperheight=228mm,paperwidth=160mm,margin=25mm]{geometry}
\allowdisplaybreaks

\theoremstyle{plain}
\newtheorem{theorem}{Theorem}[section]
\newtheorem*{theorem*}{Theorem}
\newtheorem*{quillen}{Quillen's conjecture}
\newtheorem{lemma}[theorem]{Lemma}
\newtheorem{corollary}[theorem]{Corollary}
\newtheorem{proposition}[theorem]{Proposition}

\theoremstyle{definition}
\newtheorem{example}[theorem]{Example}
\newtheorem{definition}[theorem]{Definition}

\theoremstyle{remark}
\newtheorem{remark}[theorem]{Remark}

\newtheoremstyle{TheoremNum}
		{8pt}{8pt}              	
		{\itshape}                      
		{}                              
		{\bfseries}                     
		{.}                             
		{ }                             
		{\thmname{#1}\thmnote{ \bfseries #3}}
		
\theoremstyle{TheoremNum}
\newtheorem{repeatedTheorem}{Theorem}
\newtheorem{repeatedCorollary}{Corollary}

\def\co{\colon}
\def\wt{\widetilde}
\newcommand{\tq}{\mathrel{{\ensuremath{\: : \: }}}}

\def\Aut{\mathrm{Aut}}
\def\Inn{\mathrm{Inn}}
\def\Out{\mathrm{Out}}
\def\A{\mathcal{A}}
\def\B{{\mathcal B}}
\def\F{\mathcal{F}}
\def\K{\mathcal{K}}
\def\N{\mathcal{N}}
\def\S{\mathcal{S}}
\def\Z{\mathbb{Z}}
\def\SLV{\mathcal{SLV}}
\def\PSL{\mathrm{PSL}}
\def\Sz{\mathrm{Sz}}

\def\PGU{\mathrm{PGU}}
\def\normal{\triangleleft}

\def\Inndiag{\mathrm{Inndiag}}

\def\groupiso{\cong}
\def\homotequiv{\simeq}

\begin{document}

\title[Acyclic $2$-dimensional complexes and Quillen's conjecture]{Acyclic $2$-dimensional complexes and Quillen's conjecture}

\author[K.I. Piterman]{Kevin Iv\'an Piterman$^{1}$}
\author[I. Sadofschi Costa]{Iv\'an Sadofschi Costa$^{2}$}
\author[A. Viruel]{Antonio Viruel$^{3}$}

\thanks{This work was partially done at the University of M\'alaga, during a research stay of the first two authors, supported by project MTM2016-78647-P.\\
$^1$ Supported by a CONICET doctoral fellowship and grants CONICET PIP 11220170100357CO and UBACyT 20020160100081BA.\\
$^2$ Supported by a CONICET postdoctoral fellowship and grants ANPCyT PICT-2017-2806, CONICET PIP 11220170100357CO and UBACyT 20020160100081BA. \\
$^{3}$ Partially supported by Ministerio de Econom\'ia y Competitividad (Spain), grant MTM2016-78647-P (AEI/FEDER, UE, support included).}

\address{Universidad de Buenos Aires. Facultad de Ciencias Exactas y Naturales. Departamento de Matem\'atica. Buenos Aires, Argentina.}

\address{CONICET-Universidad de Buenos Aires. Instituto de Investigaciones Matem\'aticas Luis A. Santal\'o (IMAS). Buenos Aires, Argentina.}

\address{Departamento de \'Algebra, Geometr\'ia y Topolog\'ia, Universidad de M\'alaga, Campus de Teatinos, 29071 M\'alaga, Spain.}

\email{kpiterman@dm.uba.ar}
\email{isadofschi@dm.uba.ar}
\email{viruel@uma.es}

\begin{abstract}
Let $G$ be a finite group and $\A_p(G)$ be the poset of nontrivial elementary abelian $p$-subgroups of $G$.
Quillen conjectured that $O_p(G)$ is nontrivial if $\A_p(G)$ is contractible.
We prove that $O_p(G)\neq 1$ for any group $G$ admitting a $G$-invariant acyclic $p$-subgroup complex of dimension $2$.
In particular, it follows that Quillen's conjecture holds for groups of $p$-rank $3$.
We also apply this result to establish Quillen's conjecture for some particular groups not considered in the seminal  work of Aschbacher--Smith.
\end{abstract}



\maketitle

\section{Introduction}

The study of the poset $\S_p(G)$ of nontrivial $p$-subgroups of a finite group $G$ started when K.S. Brown proved that the Euler characteristic $\chi(\K(\S_p(G)))$ of its order complex is $1$ modulo the greatest power of $p$ dividing the order of $G$ \cite{Brown}.
Recall that the order complex $\K(X)$ of a poset $X$ is the simplicial complex whose simplices are the finite nonempty totally ordered subsets of $X$.
Some years later, D. Quillen  studied the homotopy properties of $\K(\S_p(G))$ \cite{Q}.
In that article, Quillen considered the subposet $\A_p(G)$ of nontrivial elementary abelian $p$-subgroups and proved that its order complex is homotopy equivalent to $\K(\S_p(G))$ \cite[Proposition 2.1]{Q}.
Quillen also proved that, if the largest normal $p$-subgroup $O_p(G)$ of $G$ is nontrivial, then $\K(\A_p(G))$ is contractible \cite[Proposition 2.4]{Q} and conjectured that the converse should hold.

In this paper we study the following version of Quillen's conjecture.
Recall that the homology of a poset is the homology of its order complex.

\begin{quillen}
If $O_p(G)=1$ then $\wt{H}_*(\A_p(G))\neq 0$.
\end{quillen}

Aschbacher and Smith's formulation relates rational acyclicity of $\mathcal{K}(\A_p(G))$ with nontriviality of $O_p(G)$ \cite{AschbacherSmith}.
Thus our integral homology version is stronger than Quillen's original statement but weaker than the Aschbacher--Smith version.

Quillen proved the conjecture for solvable groups \cite[Theorem 12.1]{Q}.
In \cite{AschbacherSmith}, M. Aschbacher and S.D. Smith made a huge progress on the study of this conjecture.
By using the classification of finite simple groups, they proved that Quillen's conjecture holds if $p>5$ and $G$ does not contain certain unitary components.
Previously, Aschbacher and Kleidman \cite{AK} had proved Quillen's conjecture for almost simple groups (i.e. finite groups $G$ such that $L\leq G\leq \Aut(L)$ for some non-abelian simple group $L$).

The main result of our paper, which depends on the classification of the finite simple groups, is the following.

\begin{repeatedTheorem}[\ref{genthm}]
If $X$ is an acyclic and $2$-dimensional $G$-invariant subcomplex of $\K(\S_p(G))$, then $O_p(G)\neq 1$.
\end{repeatedTheorem}
Recall that the action of $G$ on $\S_p(G)$ is by conjugation.
The previous result provides then a convenient tool to prove that a group verifies Quillen's conjecture.

\begin{repeatedCorollary}[\ref{mainCorollary}]
Let $G$ be a finite group.
Suppose that $\K(\S_p(G))$ admits a $2$-dimensional and $G$-invariant subcomplex homotopy equivalent to itself.
Then Quillen's conjecture holds for $G$.
\end{repeatedCorollary}

In particular, it follows that Quillen's conjecture holds for groups of $p$-rank $3$.
Recall that the \textit{$p$-rank} of $G$, usually denoted by $m_p(G)$, is the maximum possible rank of an elementary abelian $p$-subgroup of $G$.
The $p$-rank $2$ case was considered by Quillen \cite[Proposition 2.10]{Q} and is a consequence of Serre's result: an action of a finite group on a tree has a fixed point.

In Section \ref{ExamplesSection} we make an extensive use of Corollary \ref{mainCorollary} to establish Quillen's conjecture for some particular groups (of $p$-ranks $3$ and $4$) for which the hypotheses of the results of Aschbacher--Smith \cite{AschbacherSmith} do not hold.

A related conjecture, due to C. Casacuberta and W. Dicks, is that a finite group acting on a contractible $2$-complex has a fixed point \cite{CD}.
This conjecture was studied by Aschbacher and Segev in \cite{AschbacherSegev}.
Posteriorly Oliver and Segev classified the groups which admit a fixed point free action on an acyclic (finite) $2$-complex \cite{OS}.
Our proof of Theorem \ref{genthm} is built upon the results of \cite{OS}, which  depend on the classification of finite simple groups.
Theorem \ref{genthm} can also be seen as a special case of the Casacuberta--Dicks conjecture.

\textbf{Acknowledgements.} We are grateful to the anonymous referee for their suggestions which greatly improved the exposition of the paper
and in particular for simplifying the proofs in Examples \ref{Example2a} and \ref{Example2b} by indicating Proposition \ref{propDimensionBpG}.

\section{The results of Oliver and Segev}

In this section we review the results of \cite{OS} needed in the proof of Theorem \ref{genthm}.
By a \textit{$G$-complex} we mean a $G$-CW complex. Note that the order complex of a $G$-poset is always a $G$-complex.

\begin{definition}[{\cite{OS}}]
 A $G$-complex $X$ is \textit{essential} if there is no normal subgroup $1\neq N\triangleleft G$ such that for each $H\subseteq G$, the inclusion $X^{HN}\to X^H$ induces an isomorphism on integral homology.
\end{definition}

The main results of \cite{OS} are the following two theorems.

\begin{theorem}[{\cite[Theorem A]{OS}}]\label{teoA}
 For any finite group $G$, there is an essential fixed point free $2$-dimensional (finite) acyclic $G$-complex if and only if $G$ is isomorphic to one of the simple groups $\PSL_2(2^k)$ for $k\geq 2$, $\PSL_2(q)$ for $q\equiv \pm 3 \pmod 8$ and $q\geq 5$, or $\Sz(2^k)$ for odd $k\geq 3$. Furthermore, the isotropy subgroups of any such $G$-complex are all solvable.
\end{theorem}

\begin{theorem}[{\cite[Theorem B]{OS}}]\label{teoB}
 Let $G$ be any finite group, and let $X$ be any $2$-dimensional acyclic $G$-complex. Let $N$ be the subgroup generated by all normal subgroups $N'\triangleleft G$ such that $X^{N'}\neq \emptyset$. Then $X^N$ is acyclic; $X$ is essential if and only if $N=1$; and the action of $G/N$ on $X^N$ is essential.
\end{theorem}

The set of subgroups of $G$ will be denoted by $\S(G)$.

\begin{definition}[{\cite{OS}}]
 By a \textit{family} of subgroups of $G$ we mean any subset $\F\subseteq\S(G)$ which is closed under conjugation.
 A nonempty family is said to be \textit{separating} if it has the following three properties: (a) $G\notin \F$; (b) if $H'\subseteq H$ and $H\in \F$ then $H'\in \F$; (c) for any $H\triangleleft K\subseteq G$ with $K/H$ solvable, $K\in \F$ if $H\in \F$.

 For any family $\F$ of subgroups of $G$, a \textit{$(G,\F)$-complex} will mean a $G$-complex all of whose isotropy subgroups lie in $\F$. A $(G,\F)$-complex is \textit{$H$-universal} if the fixed point set of each $H\in \F$ is acyclic.
\end{definition}

\begin{lemma}[{\cite[Lemma 1.2]{OS}}]\label{lemma1.2}
 Let $X$ be any $2$-dimensional acyclic $G$-complex without fixed points. Let $\F$ be the set of subgroups $H\subseteq G$ such that $X^H\neq \emptyset$. Then $\F$ is a separating family of subgroups of $G$, and $X$ is an $H$-universal $(G,\F)$-complex.
\end{lemma}

If $G$ is not solvable, the separating family of solvable subgroups of $G$ is denoted by $\SLV$.

\begin{proposition}[{\cite[Proposition 6.4]{OS}}]\label{proposition6.4}
 Assume that $L$ is one of the simple groups $\PSL_2(q)$ or $\Sz(q)$, where $q=p^k$ and $p$ is prime ($p=2$ in the second case). Let $G\subseteq\Aut(L)$ be any subgroup containing $L$, and let $\F$ be a separating family for $G$. Then there is a $2$-dimensional acyclic $(G,\F)$-complex if and only if $G=L$, $\F=\SLV$, and $q$ is a power of $2$ or $q\equiv \pm 3 \pmod 8$. 
\end{proposition}

\begin{definition}[{\cite[Definition 2.1]{OS}}]
 For any family $\F$ of subgroups of $G$ define 
 $$i_\F(H)=\frac{1}{[N_G(H):H]}(1-\chi(\K(\F_{>H}))).$$
\end{definition}

\begin{lemma}[{\cite[Lemma 2.3]{OS}}]\label{lemma2.3}
 Fix a separating family $\F$, a finite $H$-universal $(G,\F)$-complex $X$, and a subgroup $H\subseteq G$. For each $n$, let $c_n(H)$ denote the number of orbits of $n$-cells of type $G/H$ in $X$. Then $i_\F(H)=\sum_{n\geq 0} (-1)^nc_n(H)$.
\end{lemma}

\begin{proposition}[{\cite[Tables 2,3,4]{OS}}]\label{indices}
 Let $G$ be one of the simple groups $\PSL_2(2^k)$ for $k\geq 2$, $\PSL_2(q)$ for $q\equiv \pm 3 \pmod 8$ and $q\geq 5$, or $\Sz(2^k)$ for odd $k\geq 3$. Then $i_\SLV( 1 ) = 1$.
\end{proposition}

\section[The two-dimensional case]{The two-dimensional case}

Using the results of Oliver and Segev stated in the previous section we prove the following.

\begin{theorem}\label{thmNormalStabilizer}

Every acyclic $2$-dimensional $G$-complex has an orbit with normal stabilizer.

\begin{proof}
If $X^G\neq \emptyset$ we are done.
Otherwise, $G$ acts fixed point freely on $X$.
Consider the subgroup $N$ generated by the subgroups $N'\triangleleft G$ such that $X^{N'}\neq \emptyset$.
Clearly $N$ is normal in $G$.
By Theorem \ref{teoB} $Y=X^N$ is acyclic (in particular it is nonempty) and the action of $G/N$ on $Y$ is essential and fixed point free.
By Lemma \ref{lemma1.2} $\F=\{ H\leq G/N \tq Y^H\neq \emptyset\}$ is a separating family and $Y$ is an $H$-universal $(G/N,\F)$-complex.
Thus, Theorem \ref{teoA} asserts that $G/N$ must be one of the groups $\PSL_2(2^k)$ for $k\geq 2$, $\PSL_2(q)$ for $q\equiv \pm 3 \pmod 8$ and $q\geq 5$, or $\Sz(2^k)$ for odd $k\geq 3$.
In any case, by Proposition \ref{proposition6.4} we must have $\F=\SLV$.
By Proposition \ref{indices}, $i_\SLV(1)=1$.
Finally by Lemma \ref{lemma2.3}, $Y$ must have at least one free $G/N$-orbit.
Therefore $X$ has a $G$-orbit of type $G/N$ and we are done.
\end{proof}
\end{theorem}

\begin{theorem}\label{genthm}
If $X$ is an acyclic and $2$-dimensional $G$-invariant subcomplex of $\K(\S_p(G))$, then $O_p(G)\neq 1$.
\begin{proof}
By Theorem \ref{thmNormalStabilizer} there is a simplex $\sigma=(A_0<\ldots <A_j)$ of $X$ with stabilizer $N\triangleleft G$.
Since $A_0\triangleleft N$, we deduce that $O_p(N)$ is nontrivial.
On the other hand, $N\triangleleft G$ and $O_p(N) \, \mathrm{char}\, N$ implies that $O_p(N)\triangleleft G$.
Therefore $O_p(N)\leq O_p(G)$ and $O_p(G)$ is thus nontrivial.
\end{proof}

\end{theorem}

From Theorem \ref{genthm} we deduce:

\begin{corollary}\label{mainCorollary}
Let $G$ be a finite group.
Suppose that $\K(\S_p(G))$ admits a $2$-dimensional and $G$-invariant subcomplex homotopy equivalent to itself.
Then Quillen's conjecture holds for $G$.
\end{corollary}

Since the $p$-rank of $G$ is equal to $\dim \K(\A_p(G))+1$ we obtain:

\begin{corollary}\label{prank3}
Let $G$ be a finite group of $p$-rank $3$. If $\wt{H}_*(\A_p(G))=0$ then $O_p(G)\neq 1$.
\end{corollary}

We now apply Corollary \ref{mainCorollary} to obtain results for some related $p$-subgroup complexes.
Recall that a $p$-subgroup $Q\leq G$ is \textit{radical} if $Q = O_p(N_G(Q))$.
The \textit{Bouc poset} $\B_p(G)$ is the poset of nontrivial radical $p$-subgroups of $G$.
It is well-known that $\K(\B_p(G))$ is homotopy equivalent to $\K(\S_p(G))$ \cite{Bouc}.
Then by Corollary \ref{mainCorollary} we have

\begin{corollary}\label{coroBpG}
Let $G$ be a finite group such that $\B_p(G)$ has height $2$. If $\wt{H}_*(\B_p(G))=0$ then $O_p(G)\neq 1$.
\end{corollary}

We say that a poset $X$ is a \textit{reduced lattice} if it is obtained from a finite lattice by removing its minimum and maximum.
If $X$ is a reduced lattice, $\mathfrak{i}(X)$ denotes the subposet of $X$ given by the elements which can be written as the infimum of a set of maximal elements of $X$.
It is a general fact that the order complex of $\mathfrak{i}(X)$ is homotopy equivalent to the order complex of $X$ for any reduced lattice $X$ \cite[Section 9.1]{BarmakBook}.
Hence by Corollary \ref{mainCorollary}, we have

\begin{corollary}\label{coroiApG}
Let $G$ be a finite group. If either $\mathfrak{i}(\S_p(G))$ or $\mathfrak{i}(\A_p(G))$ has height $2$, then $G$ satisfies Quillen's conjecture.
\end{corollary}

For a detailed account of the relations between the different $p$-subgroup complexes, see \cite{S}.

\section{Some examples}\label{ExamplesSection}

In this section we apply the corollaries of Theorem \ref{genthm} to establish Quillen's conjecture for some groups constructed so that the hypotheses of the results of \cite{AschbacherSmith} are not satisfied.
The main result of \cite{AschbacherSmith} is the following.
\begin{theorem}[Aschbacher--Smith {\cite[Main Theorem]{AschbacherSmith}}]\label{ThmAschbacherSmith}
Let $G$ be a finite group and $p > 5$ a prime number.
Assume that whenever $G$ has a unitary component $U_n(q)$ with $q\equiv -1 \mod p$ and $q$ odd,
then the Quillen dimension property at $p$ holds for all $p$-extensions of $U_m(q^{p^e})$ with $m\leq n$ and $e\in \Z$.
Then $G$ satisfies Quillen's conjecture.
\end{theorem}
Recall that a group $H$ satisfies the \textit{Quillen dimension property at $p$} if $\tilde{H}_{m_p(H)-1}(\A_p(H))\neq 0$.
The presence of simple components of $G$ isomorphic to $L_2(2^3)$ or $U_3(2^3)$ (in the $p=3$ case) and $\Sz(2^5)$ (in the $p=5$ case) is an obstruction to extending Theorem \ref{ThmAschbacherSmith} to $p = 3$ and $p = 5$.
The case $p = 2$ is not considered in \cite{AschbacherSmith} and would require a much more detailed analysis.
One of the first steps in the proof of Theorem \ref{ThmAschbacherSmith} is the reduction to the case $O_{p'}(G) = 1$ (see \cite[Proposition 1.6]{AschbacherSmith}).
To do this, \cite[Theorems 2.3 and 2.4]{AschbacherSmith} are needed and these theorems make a strong use of the hypothesis $p>5$.
Concretely, it is not possible to apply \cite[Theorem 2.3]{AschbacherSmith} if a component of $C_G(O_{p'}(G))$ is isomorphic to $L_2(2^3)$, $U_3(2^3)$ (if $p=3$) or $\Sz(2^5)$ (if $p=5$).

Before presenting the examples for $p = 3$ and $p = 5$, we give some motivation.
Most of the groups $G$ in these examples satisfy the following conditions.
First, $O_{p'}(G) \neq 1$ and $C_G(O_{p'}(G))$ contains a component isomorphic to $U_3(2^3)$ if $p = 3$ and to $\Sz(2^5)$ if $p = 5$.
Thus, we cannot find nontrivial homology for $\A_p(G)$ in the same way it is done in the proof of \cite[Proposition 1.6]{AschbacherSmith} since we are not able to invoke \cite[Theorems 2.3 and 2.4]{AschbacherSmith}.
Secondly, since there is an inclusion $\tilde{H}_*(\A_p(G/O_{p'}(G));\mathbb{Q})\hookrightarrow \tilde{H}_*(\A_p(G);\mathbb{Q})$ (see \cite[Lemma 0.12]{AschbacherSmith}), we require $O_p(G/O_{p'}(G)) \neq 1$ so that $\tilde{H}_*(\A_p(G/O_{p'}(G))) = 0$.
Finally, we require $O_p(G) = 1$.

The groups presented in Examples \ref{example3Prank3} and \ref{example5Prank3} have $p$-rank $3$.
The groups presented in Examples \ref{example3Prank4} and \ref{example5Prank4} have $p$-rank $4$ and are constructed in the following way.
We take a direct product of a group $N$, consisting of one or more copies of a particular simple $p'$-group, by a group $K$ consisting of one or more copies of $L = U_3(2^3)$ if $p = 3$ or $L = \Sz(2^5)$ if $p = 5$.
Then we take two cyclic $p$-groups $A$ and $B$ and we let them act on the direct product $N\times K$ as follows.
We take a faithful action of $A\times B$ on $N$, and we choose a representation $A\times B\to \Aut(K)$ such that $O_p(K\rtimes (A\times B))\groupiso O_p( C_A(K))\neq 1$.
The group $G = (N\times K)\rtimes (A\times B)$ satisfies the conditions $O_p(G) = 1$, $O_{p'}(G) = N \neq 1$, $C_G(N) = K$ and $O_p(G/N) = O_p(K\rtimes (A\times B))\neq 1$.
Moreover, since the $p$-rank of $L$ is at most $2$, we can construct $G$ to have $p$-rank $4$ by adjusting the number of copies of $L$ in $K$.

For these groups we show that $\K(\S_p(G))$ has a $2$-dimensional $G$-invariant subcomplex homotopy equivalent to itself, and thus Corollary \ref{mainCorollary} applies.

In Examples \ref{Example2a} and \ref{Example2b} we describe two groups of $2$-rank $4$ such that $\K(\S_2(G))$ admits a $2$-dimensional $G$-invariant homotopy equivalent subcomplex.

For the claims on the structure of the automorphism group of the finite groups of Lie type we refer to \cite{GLS3} and \cite{GLS4}.

\begin{lemma}\label{lemmaprank}
Let $1\to N\to G\to K\to 1$ be an extension of finite groups.
Then
$$m_p(G)= \displaystyle\max_{A\in \S} \,\, m_p(C_N(A))\hspace{-1.5pt}+\hspace{-1pt}m_p(A),$$
where $\S$ is the set of elementary abelian $p$-subgroups $1\leq A\leq G$ such that $A\cap N=1$.
In particular we have $m_p(G)\leq m_p(N)+m_p(K)$.

\begin{proof}
If $A\in \S$ we have $C_N(A)\times A\groupiso C_N(A)A$ and hence $m_p(C_N(A))+m_p(A)\leq m_p(C_N(A)A)\leq m_p(G)$.
Taking maximum over $A\in \S$ gives the lower bound for $m_p(G)$. 
We now prove the other inequality.
Let $E$ be an elementary abelian $p$-subgroup of $G$ and write $E = (E\cap N)A$ for some complement $A$ of $E\cap N$ in $E$.
Then $m_p(E\cap N)\leq m_p(C_N(A))$ and $A\in \S$.
Now $m_p(E) = m_p(E\cap N) + m_p(A)\leq m_p(C_N(A))+m_p(A)$, giving the upper bound for $m_p(G)$.
For the last claim note that $C_N(A)\leq N$ and $m_p(A)\leq m_p(K)$ by the isomorphism theorems.
\end{proof}
\end{lemma}

The following lemma will be used to obtain proper subcomplexes of $\K(\A_p(G))$ without changing the homotopy type.
We write $X\simeq Y$ if the order complexes $\K(X)$ and $\K(Y)$ are homotopy equivalent.

\begin{lemma}\label{lemmaRetract}
Let $G$ be a finite group and let $H\leq G$.
In addition, suppose that $O_p(C_H(E))\neq 1$ for each $E\in\A_p(G)$ with $E \cap H = 1$.
Then $\A_p(G)\homotequiv \A_p(H)$.

\begin{proof}
Consider the subposet $\N = \{E\in\A_p(G) \tq E\cap H\neq 1\}$.
We have order preserving maps $r\co \N\to\A_p(H)$ and $i\co \A_p(H)\hookrightarrow\N$ given by $r(E) = E\cap H$ and $i(E) = E$ such that $ir(E)\leq E$ and $ri(E) = E$.
Therefore $\N\simeq \A_p(H)$.

Let $\S=\{E\in \A_p(G)\tq E\cap H=1 \}$ be the complement of $\N$ in $\A_p(G)$.
For any $E\in \S$ consider $\A_p(G)_{>E}\cap \N = \{A\in \N \tq A >E\}$.
It is easy to see that $r\co \A_p(G)_{>E}\cap \N \to \A_p(C_H(E))$ defined by $r(B) = B\cap H$ is a homotopy equivalence with inverse $i(B) = BE$.
Then $\A_p(G)_{>E}\cap \N\simeq \A_p(C_H(E))$ is contractible since $O_p(C_H(E))\neq 1$.

Now take a linear extension $E_1,\ldots ,E_r$ of $\S$ (i.e. ennumerate the elements of $\S$ so that $E_i\leq E_j$ implies $i\leq j$) and let $X^i = \N \cup \{E_1,\ldots,E_i\}$.
Note that $X^i = X^{i-1}\cup \{E_i\}$ and by the linear extension $X^i_{>E_i} = \A_p(G)_{> E_i}\cap \N$, which is contractible.
Now $X^i_{\geq E_i}$ is a cone over $X^i_{>E_i}$ with vertex $E_i$.
Therefore $X^{i-1} \hookrightarrow X^{i}$ is a homotopy equivalence for each $1\leq i\leq r$.
In consequence,
\begin{equation*}
\A_p(G) =X^r \simeq X^0 = \N\simeq \A_p(H). \qedhere 
\end{equation*}

\end{proof}
\end{lemma}

\begin{remark}
In the above result it can be shown that if $H\triangleleft G$ then the homotopy equivalence is $G$-equivariant.
\end{remark}

\begin{example}\label{example3Prank3}
Let $p  =3$ and let $ L = L_2(2^3) \times L_2(2^3)\times L_2(2^3)$.
Let $A$ be a cyclic group of order $3$ acting on $L$ by permuting the copies of $L_2(2^3)$.
Take $G = L\rtimes A$.
Since $m_3(L_2(2^3)) = 1$ and $C_L(A) \groupiso L_2(2^3)$, we see that $m_3(G) = 3$.
By Corollary \ref{prank3}, $G$ satisfies Quillen's conjecture.

\end{example}

\begin{example}\label{example3Prank4}
Let $p = 3$, $N = \Sz(2^3)\times \Sz(2^3)\times \Sz(2^3)$ and $U = U_3(2^3)$.
Let $A = \langle a\rangle$ and $B = \langle b\rangle$ be cyclic groups of order $3$.
We construct a semidirect product $G = (N\times U)\rtimes (A\times B)$.
To do this we need to define a map $A\times B\to \Aut(N\times U) = \Aut(N)\times \Aut(U)$.

Choose a field automorphism $\phi\in\Aut(U_3(2^3))$ of order $3$. 
By the properties of the $p$-group actions, there exists an inner automorphism $x\in \Inn(U_3(2^3))$ of order $3$ commuting with $\phi$.
Then $A\times B\to \Aut(U_3(2^3))$ is given by $a\mapsto x$ and $b\mapsto \phi$.
Choose a field automorphism $\psi\in\Aut(\Sz(2^3))$ of order $3$.
Let $A$ act on each coordinate of $N$ as $\psi$ and let $B$ act on $N$ by permuting its coordinates.
This gives rise to a well defined map $A\times B\to \Aut(N)$.

The $3$-rank of $G$ is $m_3(G) = m_3(U_3(2^3)AB)$.
We can take an elementary abelian subgroup $E\leq C_U(\phi)$ of order $9$ containing $x$ since $C_U(\phi) \groupiso \PGU_3(2)\groupiso ((C_3\times C_3) \rtimes Q_8)\rtimes C_3$ by \cite[Chapter 4, Lemma 3.10]{GLS4} and $\A_3(\PGU_3(2))$ is connected of height $1$.
Then $EAB$ is an elementary abelian subgroup of order $3^4$.
Hence, $m_3(UAB)\geq 4$.
Since $m_3(U_3(2^3)) = 2$ and $m_3(AB) = 2$, by Lemma \ref{lemmaprank} we have $m_3(G)=4$.

By Corollary \ref{mainCorollary}, to show that Quillen's conjecture holds for $G$ and $p = 3$ it is enough to find a $2$-dimensional $G$-invariant subcomplex $X$ of $\K(\S_3(G))$ homotopy equivalent to $\K(\S_3(G))$ (or, equivalently, to $\K(\A_3(G))$).

Let $H = (N\times U) \rtimes A$.
Note that $H\triangleleft G$ and $m_3(H) = 3$.
Therefore, $\K(\A_3(H))$ is a $2$-dimensional $G$-invariant subcomplex of $\K(\A_3(G))$.
Now the plan is to use Lemma \ref{lemmaRetract} to show that $\A_3(H)\simeq \A_3(G)$.
Let $E \in \A_3(G)$ be such that $E\cap H = 1$.
Then $E\groupiso EH/H \leq B \groupiso C_3$ and hence, $E$ is cyclic generated by some element $e\in E$.
Write $e = nua^ib^j$ with $n\in N$, $u\in U$ and $i,j\in \{0,1,2\}$.
Note that $j\neq 0$ since $E\cap H = 1$.
If $v\in U$, then
$$v^{e} = v^{nua^ib^j} = (v^{ua^i})^{b^j}.$$
Since $j\neq 0$ and $e$ induces an automorphism of $U$ of order $3$ in $\Inn(U) \phi^j $, by \cite[Proposition 4.9.1]{GLS3} and the definition of field automorphisms \cite[Definition 2.5.13]{GLS3}, $e$ is $\Inndiag(U)$-conjugate to $\phi^j$ and acts as a field automorphism on $U$.
In particular, $C_U(E) = C_U(e) \groupiso C_U(\phi^j) = C_U(\phi)$.
Note that $O_3(C_U(E)) \groupiso O_3(C_U(\phi)) \groupiso C_3\times C_3\neq 1$.
Since $C_U(E)\triangleleft C_H(E)$ and $O_3(C_U(E))\neq 1$, we conclude that $O_3(C_H(E))\neq 1$.
By Lemma \ref{lemmaRetract}, $\A_3(G)\homotequiv \A_3(H)$, which is $2$-dimensional and $G$-invariant.
In conclusion, the subcomplex $\K(\A_3(H))$ satisfies the hypothesis of Corollary \ref{mainCorollary} and therefore, Quillen's conjecture holds for $G$.

Note that $O_3(G) = 1$, $O_{3'}(G) = N$, $C_G(O_{3'}(G)) = U_3(2^3)$ and $O_3(G / O_{3'}(G)) = O_3(U_3(2^3)AB) = \langle ax^{-1}\rangle \groupiso C_3$.
\end{example}

\begin{example}\label{example5Prank3}
Let $p=5$.
Let $r$ be a prime number such that $r\equiv 2$ or $3\mod 5$ and let $q = r^{5^n}$ with $n\geq 2$.
Let $N$ be one of the simple groups $L_2(q)$, $G_2(q)$, ${}^3D_4(q^3)$ or ${}^2G_2(3^{5^n})$  and let $A=\langle a \rangle$ be a cyclic group of order ${5^n}$.
Note that $5\nmid |N|$.
Let $a$ act on $N$ as a field automorphism of order $5^n$.
Choose a field automorphism $\phi\in\Aut(\Sz(2^5))$  of order $5$ and let $A$ act on $\Sz(2^5)\times \Sz(2^5)$ as $\phi\times\phi$.
Now consider the semidirect product $G = (N\times \Sz(2^5)\times \Sz(2^5))\rtimes A$ defined by this action.

Since the Sylow $5$-subgroups of $\Sz(2^5)$ are cyclic of order $25$, by Lemma \ref{lemmaprank} we have that $m_5(G) = 3$.
By Corollary \ref{prank3}, Quillen's conjecture holds for $G$.

Moreover, $O_5(G) = 1$, $O_{5'}(G) = N$, $C_G(O_{5'}(G)) = \Sz(2^5)^2$ and $O_5(G/O_{5'}(G)) = C_A(\Sz(2^5)^2) = \langle a^5\rangle \neq 1$.
\end{example}

\begin{example}\label{example5Prank4}
Let $p = 5$ and let $N = L^5$, where $L$ is one of the simple $5'$-groups of the previous example.
Let $A =\langle a \rangle \groupiso C_{5^n}$ and $B=\langle b \rangle \groupiso C_5$.
Let $G = (N \times \Sz(2^5)^2)\rtimes (A\times B)$, where $a$ acts on each copy of $L$ as a field automorphism of order $5^n$ and trivially on $\Sz(2^5)^2$, and $b$ permutes the copies of $L$ and acts as a field automorphism of order $5$ on each copy of $\Sz(2^5)$.

To compute the $5$-rank of $G$ we use Lemma \ref{lemmaprank}:
\begin{align*}
m_5(G) 	  &= m_5(\Sz(2^5)^2\rtimes (A\times B)) \\
	  &=  m_5(A\times (\Sz(2^5)^2\rtimes B)) \\
	  &= m_5(A) + m_5(\Sz(2^5)^2\rtimes B) \\
	  &= 1 + 3\\
	  &= 4.
\end{align*}
Now the aim is to apply Corollary \ref{mainCorollary} on $G$ by finding a $2$-dimensional $G$-invariant homotopy equivalent subcomplex $X$ of $\K(\S_5(G))$.

Let $H = (N\times \Sz(2^5)^2)\rtimes A = NA\times \Sz(2^5)^2$.
Note that $H\triangleleft G$ and $m_5(H) = 3$. Hence $\K(\A_5(H))$ is $2$-dimensional and $G$-invariant.
We will show that $\A_5(H)\homotequiv \A_5(G)$ by applying Lemma \ref{lemmaRetract}.

Let $E\in \A_5(G)$ be such that $E\cap H = 1$.
Then $E$ is cyclic generated by an element $e$ of order $5$ and $e = lsa^ib^j$ with $l\in N$, $s\in\Sz(2^5)^2$, $0\leq i\leq 5^n-1$ and $j\in\{1,2,3,4\}$.
Thus $E$ acts by field automorphisms on each copy of the Suzuki group and $e$ is $\Inndiag(\Sz(2^5))$-conjugate to the field automorphism induced by $b^j$ on $\Sz(2^5)$ (see \cite[Proposition 4.9.1]{GLS3} and Example \ref{example3Prank4}).
Hence, $C_H(E) = C_{NA}(E) \times C_{\Sz(2^5)^2}(E)$.
Note that  $C_{\Sz(2^5)^2}(E)\normal C_H(E)$ and $C_{\Sz(2^5)^2}(E) \groupiso C_{\Sz(2^5)}(E)^2\groupiso (C_5\rtimes C_4)^2$ has a nontrivial normal $5$-subgroup.
Therefore $\A_5(G)\homotequiv \A_5(H)$ by Lemma \ref{lemmaRetract} and Quillen's conjecture holds for $G$ by Corollary \ref{mainCorollary} applied to the subcomplex $\K(\A_5(H))$.

Note that $O_{5'}(G) = N$ and $C_G(O_{5'}(G)) = \Sz(2^5)^2$. On the other hand, $O_5(G) = 1$ and $O_5(G/O_{5'}(G)) = A\neq 1$.
\end{example}

We conclude with two examples of groups satisfying Quillen's conjecture for $p = 2$.
We say that a finite group $G$ has the \textit{trivial intersection property at $p$} if any two different Sylow $p$-subgroups of $G$ have trivial intersection.

\begin{proposition}\label{propDimensionBpG}
Let $L_1$ and $L_2$ be two finite groups with the trivial intersection property at $p$.
Let $L = L_1\times L_2$ and take an extension $G$ of $L$ such that $|G:L| = p$.
Then $\mathfrak{i}(\S_p(G))$ and $\B_p(G)$ are at most $2$-dimensional.
If in addition the Sylow $p$-subgroups of $L_1$ and $L_2$ have abelian $\Omega_1$, then $\mathfrak{i}(\A_p(G))$ is at most $2$-dimensional.
\begin{proof}
The elements of $\mathfrak{i}(\S_p(L))$ are of the form $P_1\times P_2$, $1\times P_2$ or $P_1\times 1$, where  $P_i\leq L_i$ are Sylow $p$-subgroups.
Hence, $\mathfrak{i}(\S_p(L))$ is $1$-dimensional.

Now suppose that $Q_0 < Q_1 < \ldots < Q_n$ is a chain in $\mathfrak{i}(\S_p(G))$.
Then
$$Q_0\cap L\leq Q_1\cap L \leq \ldots\leq Q_n\cap L$$
is a chain in $\mathfrak{i}(\S_p(L))$.
We claim that there is at most one index $i$ such that $Q_i\cap L= Q_{i+1}\cap L$.
To see this note that $$|Q_{j}:Q_j\cap L|=\begin{cases} 1 & \text{ if } Q_j\subseteq L \\ p & \text{ if } Q_j\not\subseteq L\end{cases}.$$
We have $|Q_{i+1}:Q_i|\cdot |Q_i:Q_i\cap L| = |Q_{i+1}:Q_{i+1}\cap L|\cdot|Q_{i+1}\cap L : Q_{i}\cap L|$.
Then if $Q_i\cap L = Q_{i+1}\cap L$, since $|Q_{i+1}:Q_i|\geq p$ we must have $|Q_{i}:Q_i\cap L|=1$ and $|Q_{i+1}:Q_{i+1}\cap L|=p$. Then $i=\max\{j\tq Q_j\subseteq L\}$.

From this we conclude that $\dim \mathfrak{i}(\S_p(G)) \leq 1 + \dim \mathfrak{i}(\S_p(L)) = 2$.
It is well-known that $\B_p(G)$ is a subposet of $\mathfrak{i}(\S_p(G))$ (i.e. every radical $p$-subgroup is an intersection of Sylow $p$-subgroups). Then $\B_p(G)$ is   at most $2$-dimensional also.
The same proof can be easily adapted to prove that, if the Sylow $p$-subgroups of $L_1$ and $L_2$ have abelian $\Omega_1$, $\mathfrak{i}(\A_p(G))$  is at most $2$-dimensional.
\end{proof}
\end{proposition}

In the following examples we use the fact that the groups $A_5$ and $U_3(2^2)$ have the trivial intersection property at $2$ and that $\Omega_1(P)$ is abelian for $P$ a Sylow $2$-subgroup of either $A_5$ or $U_3(2^2)$.

\begin{example}\label{Example2a}
Let $G$ be the group extension $(A_5\times A_5)\rtimes C_2$ where the generator of $C_2$ acts on each coordinate as conjugation by the transposition $(1 \,2)$.
Since $m_2(A_5) = 2 = m_2(\Aut(A_5))$, by Lemma \ref{lemmaprank}, $G$ has $2$-rank $4$.
By Proposition \ref{propDimensionBpG}, $\mathfrak{i}(\A_2(G))$, $\mathfrak{i}(\S_2(G))$ and $\B_2(G)$ are $2$-dimensional and then Quillen's conjecture holds for $G$ since Corollaries  \ref{coroBpG} and \ref{coroiApG} apply.
\end{example}

\begin{example}\label{Example2b}
Let $G=(U_3(2^2)\times A_5)\rtimes C_2$ be the semidirect product constructed in the following way.
Let $H = U_3(2^2)\times A_5$.
Then $\Out(H) \groupiso \Aut(U_3(2^2))/\Inn(U_3(2^2)) \times \Aut(A_5)/\Inn(A_5) \groupiso C_4\times C_2$.
Take $t\in \Out(H)$ to be the involution which acts nontrivially on both factors.
Therefore $G = H \rtimes \langle t\rangle$.
Since $m_2(U_3(2^2)) = 2 = m_2(A_5) = m_2(\Aut(A_5))$ and $m_2(\Aut(U_3(2^2))) = 3$, by Lemma \ref{lemmaprank} $G$ has $2$-rank $4$.
Just as before, Quillen's conjecture holds for $G$.
\end{example}

\bibliographystyle{alpha}
\bibliography{references}
\end{document}